\documentclass[a4paper,10pt]{amsart}

\usepackage{amsmath,amssymb,amsthm,mathtools}
\usepackage{booktabs}
\usepackage[hidelinks]{hyperref}

\newtheorem{theorem}{Theorem}[section]
\newtheorem{lemma}[theorem]{Lemma}
\newtheorem{proposition}[theorem]{Proposition}
\newtheorem{corollary}[theorem]{Corollary}
\theoremstyle{remark}
\newtheorem{remark}[theorem]{Remark}

\title[Power projections, Gram structures, and finite-field lifting]{Power projections, Gram structures, and finite-field lifting of Butson Hadamard matrices}

\author{Wojciech Bruzda}
\address{W.~B.: Center for Quantum-Enabled Computing, Center for Theoretical Physics, Polish Academy of Sciences, Aleja Lotników 32/46, 02-668 Warsaw, Poland}
\email{w.bruzda@cft.edu.pl}
\author{Domonkos Czifra}
\address{D.~C.: HUN-REN Alfréd Rényi Institute of Mathematics, Reáltanoda u. 13--15, \mbox{Budapest} 1053, Hungary}
\email{czifra.domonkos@renyi.hu}
\author{Ilona Nagy}
\address{I.~N.: HUN-REN Alfréd Rényi Institute of Mathematics, Reáltanoda u. 13--15, \mbox{Budapest} 1053, Hungary}
\curraddr{Department of Analysis and Operations Research, Institute of Mathematics, Budapest University
of Technology and Economics, Műegyetem rkp. 3, Budapest 1111, Hungary}
\email{nagyi@math.bme.hu}
\author{Ferenc Sz\"oll\H{o}si}
\address{F.~Sz.: Interdisciplinary Department of Science and Engineering, Shimane University, 1060
Nishikawatsu-cho, Matsue, Shimane 690-8504, Japan}
\email{szollosi@riko.shimane-u.ac.jp}
\date{November 18, 2026}

\begin{document}
\begin{abstract}
In this paper we investigate the algebraic properties of Hadamard powers of Butson Hadamard matrices $H\in \mathrm{BH}(n,3p)$, and show that for certain parameter values the cube root projection matrix $H^{\circ p}$ has the same structure as the adjacency matrix of a strongly regular graph. Then, by using a simple eigenvalue analysis we show, among other results, the nonexistence of $\mathrm{BH}(17,15)$ matrices.

For other parameter values, a more in-depth study and in particular, a complete classification of all possible projection matrices is necessary. We carry out this task for the case of $\mathrm{BH}(16,15)$ matrices, and using finite field techniques we show that none of these projections can be lifted up to form a Butson matrix. In particular, we show the nonexistence of $\mathrm{BH}(16,15)$ matrices for the first time.
\end{abstract}
\maketitle
\section{Introduction}
Let $\mathbb{C}^n$ denote the usual complex $n$-dimensional Hilbert space equipped with the standard inner product $\left\langle\cdot,\cdot\right\rangle$. A Butson Hadamard matrix $H\in\mathrm{BH}(n,q)$ is an $n\times n$ complex matrix with entries in the multiplicative group of $q$th roots of unity $\Omega_q:=\{z\in\mathbb{C}\colon z^q=1\}$, satisfying $HH^\ast=nI_n$. In other words, the rows and columns of $H$ are pairwise complex orthogonal. Equivalently, Butson matrices attain the Hadamard determinant bound \cite{NP23} among matrices over $\Omega_q$. This class of matrices was introduced by Butson in order to generalize the concept of real Hadamard matrices \cite{Butson1962}. Two matrices over $\Omega_q$ are called equivalent if one can be obtained from the other by a sequence of row and column permutations and multiplication of rows and columns by elements of $\Omega_q$. Notable examples include the Fourier matrices and their Kronecker products.

Butson matrices appear in many areas of mathematics and physics, including algebraic combinatorics \cite{MI}, operator algebras \cite{BN}, and quantum information theory \cite{TZ}. In theoretical physics, they are often rescaled by a factor of $1/\sqrt{n}$ to become unitary. They play a crucial role in quantum permutation algebras~\cite{BBS09}, unitary error bases \cite{KR02}, dense coding schemes~\cite{W01} and constructions of quantum error-correcting codes~\cite{E23}. Furthermore, a rescaled complex Hadamard matrix can be viewed as the transition matrix between two mutually unbiased orthonormal bases~\cite{DEBZ10}, and as a link between a pair of orthogonal maximal abelian $\ast$-subalgebras~\cite{P83}. Finally, imposing additional structural conditions in the form of multiunitarity leads to applications to absolutely maximally entangled states~\cite{BZ24,RBBMLZ22}.

Determining the pairs $(n,q)$ for which the class of $\mathrm{BH}(n,q)$ is nonempty remains a challenging open problem. In the simplest case of $q=2$ the famous Hadamard conjecture predicts the existence of a $\mathrm{BH}(n,2)$ matrix for every doubly even $n$, but this is far from being resolved despite continuous efforts \cite{review}. An intriguing application of Butson matrices is that in certain cases they can be lifted to higher order real Hadamard matrices \cite{CCD, EOC19, T01}. Conversely, if $k$ is the order of a real Hadamard matrix and $p$ is a large enough prime number, then there exist $\mathrm{BH}(kp,p)$ matrices \cite{DDL94}.

Considering nonexistence, important tools include divisibility and self-conjugacy conditions~\cite{B88,DL84,W00}, the structure of vanishing sums of roots of unity~\cite{LamLeung}, and exhaustive classification~\cite{LampioEtAl}. 

Several recent papers have advanced these questions in various directions. \mbox{Czifra}, Matolcsi, and Sz\"oll\H{o}si introduced a new nonexistence criterion based on positive definite functions on finite groups, and constructed new Butson matrices using cyclotomic cosets \cite{CzifraEtAl}. Leung, Chue and Zhao applied vanishing-sum techniques to cyclic Butson matrices \cite{LCZ25}, and Schmidt with coauthors classified matrices which are invariant under certain abelian groups \cite{DDS26, SWX21}. Finally, generalized bent vectors~\cite{AEKO26}, and switching operations for Butson Hadamard matrices have also been introduced recently \cite{CES26}.

Let $\circ$ denote the Hadamard product (i.e., entrywise product) of matrices. In this work we study the algebraic properties of Hadamard powers $H^{\circ p}$ of $\mathrm{BH}(n,3p)$ matrices, mainly focusing on the case $3p=15$ which is the smallest odd root order not being a prime power. In Section~\ref{sec2} we prove in Theorem~\ref{T29} the nonexistence of $\mathrm{BH}(16,15)$ matrices, and in Section~\ref{sec3} we derive the nonexistence of $\mathrm{BH}(17,15)$ matrices in Corollary~\ref{cor}. The common starting point of the proof of these results is taking the $p$th Hadamard power of a hypothetical Butson matrix $H$ which projects its entries from $\Omega_{3p}$ to $\Omega_3$. Although the projected cube root matrix $H^{\circ p}$ can not be a Butson matrix for obvious divisibility reasons, its Gram matrix remains highly structured due to the unique vanishing-sum decomposition \cite{LamLeung} of the quotient rows of $H$. For $n=17$, the Gram matrix of $H^{\circ p}$ can be related to the adjacency matrix of a nonexistent strongly regular graph, while for $n=16$ a single projection matrix was found up to equivalence (see Proposition~\ref{prop26}), which could not be lifted back to $\Omega_{15}$ to form a genuine $\mathrm{BH}(16,15)$ matrix.

The online Butson Home catalogue\footnote{\url{ https://wiki.aalto.fi/spaces/Butson/pages/120482304/Matrices+up+to+monomial+equivalence}} records up-to-date information regarding the existence of $\mathrm{BH}(n,q)$ matrices for small parameter values $(n,q)$.
\section{The nonexistence of BH(16,15) matrices}\label{sec2}
The goal of this section is to prove Theorem~\ref{T29}, which states that $\mathrm{BH}(16,15)$ matrices do not exist. First we recall the following folklore fact from linear algebra.
\begin{lemma}[The generalized matrix determinant lemma]\label{mdl}
Let $U$ and $V$ be complex matrices of order $n\times m$ and $m\times n$, respectively. Then $\mathrm{det}(I_n+UV)=\mathrm{det}(I_m+VU)$.
\end{lemma}

Next, consider $H\in\mathrm{BH}(n,q)$, and let $e\in\Omega_{q}^{n}$ be the constant $1$ column vector of appropriate length, and for $i\in\{1,\dots,n\}$ let $H_i$ denote the $i$th row of $H$. We say that $H$ is in normal form, if $H_1=e^T$, and $\prod_{k=1}^{n}H_{ik}=1$ for every $i$. 

It is easy to see that every $H\in\mathrm{BH}(16,15)$ matrix is equivalent to one in normal form. Indeed, one can multiply the columns of $H$ by appropriate roots of unity, and then row $H_i$ by $\prod_{k=1}^{16}H_{ik}^{-1}$ to obtain a matrix in normal form. Throughout this section we assume that $H$ is in normal form.

Let $\omega\in\Omega_{15}$ be a primitive $3$rd, and $\zeta\in\Omega_{15}$ be a primitive $5$th root of unity, respectively. We consider $A:=H^{\circ 5}$, which is a cube root matrix in normal form.
\begin{lemma}\label{l2}
Let $H\in\mathrm{BH}(16,15)$ be in normal form, and let $A:=H^{\circ 5}$. Then, the Gram matrix $AA^\ast$ has off-diagonal entries $-5$ or $10$.
\end{lemma}
\begin{proof}
Since the rows of $H$ are pairwise complex orthogonal, the quotient row $H_i\circ \overline{H_j}$ ($i<j\in\{1,2,\dots,16\}$) is necessarily a permutation of
\[[h_1,h_1\zeta,h_1\zeta^2,h_1\zeta^3,h_1\zeta^4,h_2,h_2\zeta,h_2\zeta^2,h_2\zeta^3,h_2\zeta^4,h_3,h_3\omega,h_3\omega^2,h_4,h_4\omega,h_4\omega^2],\]
with some $h_1,h_2,h_3,h_4\in \Omega_{15}$, i.e., the inner product $\left\langle H_i,H_j\right\rangle$ decomposes as a vanishing sum of two $5$-cycles and two $3$-cycles \cite{LamLeung}. Consequently, $(H_i\circ \overline{H_j})^{\circ 5}=A_i\circ \overline{A_j}$ is a permutation of
\begin{equation}\label{rab}
R(\alpha,\beta):=[\alpha,\alpha,\alpha,\alpha,\alpha,\beta,\beta,\beta,\beta,\beta,1,\omega,\omega^2,1,\omega,\omega^2], \alpha,\beta\in\Omega_3,\end{equation}
and $\left\langle A_i,A_j\right\rangle = \sum_{k=1}^{16}(A_i\circ \overline{A_j})_k=5\alpha+5\beta$.

Since $A$ is in normal form, $1=\prod_{k=1}^{16} A_{ik}\cdot \prod_{k=1}^{16} \overline{A_{jk}}=\prod_{k=1}^{16} (A_i\circ \overline{A_j})_k=\alpha^5\beta^5$, equivalently $\beta=\overline{\alpha}$. In particular, $
\alpha=\beta=1$ or $\alpha=\omega$ and $\beta=\omega^2$. Hence $\left\langle A_i,A_j\right\rangle = 5\alpha+5\beta=5(\alpha+\overline{\alpha})=10 \Re[\alpha]\in \{-5,10\}$ as claimed.
\end{proof}
For $k\geq 1$ we denote by $J_k$ the constant $1$ matrix of order $k$, and omit the subscript if its size is clear from context.
\begin{lemma}\label{l2333}
We have $AA^\ast=6I-5J+15\operatorname{diag}(J_a,J_b,J_c)$ for some $(a,b,c)\in\mathbb{Z}^3$ with $a\geq b\geq c\geq 0$, $a+b+c=16$, up to row permutations. Consequently $\operatorname{det}(AA^\ast)=3\cdot 6^{14}\cdot (332+25(ab+ac+bc))$.
\end{lemma}
We note that here blocks of size $0$ are also allowed, using the convention that $\mathrm{diag}(J_a,J_b,J_0):=\mathrm{diag}(J_a,J_b)$, etc.
\begin{proof}
The matrix $AA^\ast$ is symmetric, positive semi-definite, with off-diagonal entries $10$ or $-5$ by Lemma~\ref{l2}. Therefore, it cannot contain any of the following principal submatrices having negative determinant:
\[21I_3+10J_3-15\operatorname{diag}(J_1,J_2),\quad 21I_5-5J_5,\quad 6I_7-5J_7+15\operatorname{diag}(J_4,J_1,J_1,J_1).
\]
The first restriction implies that $AA^\ast$, up to equivalence, decomposes into a block-diagonal form where blocks $6I_k+10J_k$ (of possibly different sizes $k\in\{1,2,\dots,16\}$) are on the main diagonal, and the entries are $-5$ otherwise. The second restriction implies that the number of blocks is at most $4$. Finally, the third restriction implies that the number of blocks cannot be exactly $4$. Indeed, in this case the largest block should have at least $4$ rows, hence selecting these and one each from the three other blocks would give the third forbidden submatrix. Therefore we have at most three blocks, and the claimed block decomposition holds.

To obtain the formula for the determinant, we use Lemma~\ref{mdl} where $U$ is a $\{0,1\}$-matrix of size $16\times 3$, whose columns are the indicator vectors of the blocks $J_a$, $J_b$, and $J_c$ (in particular, a constant $0$ column would correspond to $J_0$), respectively, and $V:=(15I_3-5J_3)U^T/6$. Then, we have
\begin{align*}
\operatorname{det}(AA^\ast)&=\mathrm{det}(6I+15\operatorname{diag}(J_a,J_b,J_c)-5J)=\operatorname{det}(6I+U(15I_3-5J_3)U^T)\\
&=6^{16}\operatorname{det}(I+U(15I_3-5J_3)U^T/6)=6^{16}\operatorname{det}(I_3+\frac{5}{6}(3I_3-J_3)U^TU)\\
&=6^{13}\operatorname{det}(6I+5(3I_3-J_3)\operatorname{diag}(a,b,c)).
\end{align*}
The resulting $3\times 3$ determinant can be calculated directly; replacing $a+b+c$ with its value $16$ leads to the statement.
\end{proof}
\begin{lemma}\label{l3}
We have $AA^\ast=6I-5J+15\operatorname{diag}(J_a,J_b,J_c)$ for some $(a,b,c)\in\mathbb{Z}^3$ with $a\geq b\geq c\geq 1$, $a\geq 6$, $a+b+c=16$, up to row permutations.
\end{lemma}
\begin{proof}
By Lemma~\ref{l2333}, the Gram matrix $AA^\ast$ has at most three blocks. It is easy to see that there are $30$ integer partitions of $16$ with at most three summands. However, we may eliminate most of these by observing that $\mathrm{det}(A)=x+y\omega$ for some $x,y\in\mathbb{Z}$ is an Eisenstein integer \cite[Section~4]{DAC}. Hence, the square-free part of its norm $\operatorname{det}(AA^\ast)=|x+y\omega|^2=x^2-xy+y^2$ cannot have any prime divisor $p\equiv 2\ (\operatorname{mod} 3)$. This leaves us with the following possibilities:
\begin{align*}
(a,b,c)\in{}&\{(14,1,1), (11,4,1), (11,3,2), (10,5,1), \\
{}&\phantom{\{}  (10,4,2), (8,7,1), (8,5,3), (7,7,2), (6,5,5)\}.\qedhere
\end{align*}
\end{proof}
\begin{remark}\label{rem33}
It can be shown by considering the permutations of $R(1,1)$ in formula \eqref{rab}, that the Gram matrix $AA^\ast$ can not contain the principal submatrix $6I_9+10J_9$ either. In particular, this can lead to a slightly improved variant of Lemma~\ref{l3}, where the condition $a\leq 8$ is also claimed. While we do not rely on these observations, these can be useful during the study of other pairs $(n,q)$. See also Table~\ref{tab1}.
\end{remark}
\begin{proposition}\label{prop26}
Let $H\in\mathrm{BH}(16,15)$ be in normal form, and let $A:=H^{\circ 5}$. Then, we have $AA^\ast=6I-5J+15\operatorname{diag}(J_8,J_7,J_1)$, up to row permutations. In particular, we may choose
\begin{equation}\label{AA}A=\left[\begin{array}{cccc}
1  & e^T & e^T & 1\\
\omega e & (\omega-1)I_7+J_7 & (\omega^2-1)I_7+J_7 & \omega^2e\\
e & (\omega-\omega^2) I_7 +\omega^2 J_7 & (1-\omega)I_7+\omega J_7 & \omega^2e\\
1 & \omega e^T & \omega^2 e^T & 1\end{array}\right].
\end{equation}
up to equivalence followed by a possible conjugation.
\end{proposition}
\begin{proof}
Let $x+y+z\geq 1$, and let us call an $(x+y+z)\times 16$ submatrix $X$ of $A$ of type $(x,y,z)\in\mathbb{Z}^3$, if $XX^\ast=6I-5J+15\operatorname{diag}(J_x,J_y,J_z)$. By Lemma~\ref{l3}, $A$ has an $8\times 16$ submatrix, whose Gram matrix is of type $(6,1,1)$. It is a routine task to generate these matrices up to equivalence, see Table~\ref{tab1}. Indeed, the first row can be taken as $e^T$, and then five additional rows are certain permutations of $R(1,1)$ shown in \eqref{rab}. Finally, the last two rows are permutations of $R(\omega,\omega^2)$. In the end, we find $18$ such matrices, which form only $10$ classes once complex conjugation is also considered. We store these $10$ `seed matrices' in a set $\mathcal{S}$.

Next, for each of these seeds $S\in\mathcal{S}$, we generate the rows $r\in\Omega_3^{16}$ such that $\left\langle S_i,r\right\rangle\in\{-5,10\}$ for every $i\in\{1,\dots,8\}$, and put these in a set $\mathcal{R}(S)$. If $\mathcal{R}(S)$ is empty then the seed $S$ is inextendible; otherwise we create a graph $\Gamma(S)$ whose vertex set is indexed by the elements of $\mathcal{R}(S)$. Two vertices $r_1,r_2\in\mathcal{R}(S)$ are adjacent if and only if $\left\langle r_1,r_2\right\rangle\in\{-5,10\}$. We determine the size $\omega(\Gamma(S))$ of the maximum cliques in $\Gamma(S)$, and we enumerate them in case $\omega(\Gamma(S))=8$. Using the Cliquer software \cite{patclq}, a total of $14$ matrices were found, forming a single equivalence class (where conjugate matrices are also considered equal) of type $(8,7,1)$ with representative $A$ as claimed.
\end{proof}

\begin{table}[htbp]
\centering
\caption{Summary of the computational data relevant to $A=H^{\circ 5}$}\label{tab1}
\begin{tabular}{cc}
\toprule
Type & \#\\
\midrule
$(1,0,0)$ & 1\\
$(2,0,0)$ & 1\\
$(3,0,0)$ & 6\\
$(4,0,0)$ & 10\\
$(5,0,0)$ & 10\\
$(6,0,0)$ & 13\\
$(6,1,0)$ & 24\\
$(6,1,1)$ & 18\\
$(7,0,0)$ & 7\\
$(8,0,0)$ & 3\\
$(9,0,0)$ & 0\\
\bottomrule
\end{tabular}
\qquad\qquad\qquad
\begin{tabular}{cccc}
\toprule
$S$ & $|\mathcal{R}(S)|$ & $\omega(\Gamma(S))$ & \#\\
\midrule
 1 & 0 & $-$ & $-$ \\
 2 & 6 & 1 & $-$ \\
 3 & 9 & 2 & $-$ \\
 4 & 24 & 5 & $-$ \\
 5 & 15 & 8 & 2 \\
 6 & 15 & 8 & 2 \\
 7 & 15 & 8 & 2 \\
 8 & 15 & 8 & 2 \\
 9 & 15 & 8 & 2 \\
 10 & 21 & 8 & 4 \\
\bottomrule
\\
\end{tabular}
\end{table}

Next, write $A:=[\omega^{a_{ij}}]$, $B:=[\zeta^{b_{ij}}]$, and $H:=\overline{A}\circ B$. Then $H^{\circ 5}=\overline{A}\circ \overline{A}=A$.
\begin{lemma}\label{lem4}
Let $H$ be a $\mathrm{BH}(16,15)$ matrix in normal form, and let $A$ be the matrix \eqref{AA}. Put $A_{ij}=\omega^{a_{ij}}$, and $H=\overline{A}\circ B$, $H_{ij}=\omega^{-a_{ij}}\cdot\zeta^{b_{ij}}$. For every $i<j\in\{1,2,\dots,16\}$ and for every $r\in\{0,1,2\}$, define the index set $\mathcal{I}(i,j,r):=\{k\in\{1,2,\dots,16\}\colon -(a_{ik}-a_{jk})\equiv r\ (\operatorname{mod} 3)\}$. Then, we have
\begin{align}
\sum_{k=1}^{16} b_{ik}&\equiv 0\ (\mathrm{mod}\ 5),\label{eq1}\\
\sum_{k\in \mathcal{I}(i,j,r)}(b_{ik}-b_{jk})&\equiv 0\ (\operatorname{mod} 5).\label{eqlin}
\end{align}
\end{lemma}
\begin{proof}
Since $H$ and $A$ are in normal form, so is $B$. Therefore, for every $i\in\{1,2,\dots,16\}$, we have $\prod_{k=1}^{16} \zeta^{b_{ik}}=1$, which corresponds to a linear equation of the form \eqref{eq1}.

Next, consider the inner product of the rows of $H$ for $i<j\in\{1,2,\dots,16\}$:
\[\left\langle H_i,H_j\right\rangle = \sum_{k=1}^{16}(H_i\circ \overline{H_j})_k=\sum_{r=0}^2\omega^r \sum_{k\in\mathcal{I}(i,j,r)}\zeta^{b_{ik}-b_{jk}}=0.\]
It is well-known (see \cite[Proposition 2.4]{LCW}) that $1$ and $\omega$ are linearly independent over $\mathbb{Q}(\zeta)$. Furthermore, as $\omega^2=-1-\omega$, we have (exponents are understood modulo $5$, since $\zeta^5=1$):
\begin{equation}\label{magic}\sum_{k\in\mathcal{I}(i,j,0)}\zeta^{b_{ik}-b_{jk}}=\sum_{k\in\mathcal{I}(i,j,1)}\zeta^{b_{ik}-b_{jk}}=\sum_{k\in\mathcal{I}(i,j,2)}\zeta^{b_{ik}-b_{jk}}.\end{equation}
Since the coefficients of the powers of $\zeta$ are nonnegative integers in the above sums, they must contain exactly the same terms, except for additional $N(r)$ vanishing sums of the form $1+\zeta+\zeta^2+\zeta^3+\zeta^4$, $N(r)\geq 0$, $r\in\{0,1,2\}$.

Finally, evaluate the product of the entries in the quotient rows of $B$ as follows:
\[1=\prod_{k=1}^{16}B_{ik}\cdot \prod_{k=1}^{16}\overline{B_{jk}}=\prod_{k=1}^{16}(B_{i}\circ \overline{B_j})_k=\prod_{r=0}^2\prod_{k\in\mathcal{I}(i,j,r)}\zeta^{b_{ik}-b_{jk}}.\]
Since the terms in the additional vanishing sums in \eqref{magic} have product $\zeta^{10N(r)}=1$, these do not affect the above product, hence 
\[\prod_{k\in\mathcal{I}(i,j,0)}\zeta^{b_{ik}-b_{jk}}=\prod_{k\in\mathcal{I}(i,j,1)}\zeta^{b_{ik}-b_{jk}}=\prod_{k\in\mathcal{I}(i,j,2)}\zeta^{b_{ik}-b_{jk}}\]
is a cube root of $1$. Since the only cube root in $\mathbb{Q}(\zeta)$ is $1$, \eqref{eqlin} holds.
\end{proof}
\begin{proposition}\label{prop2}
Let $H$ be a $\mathrm{BH}(16,15)$ matrix in normal form, and let $A$ be the matrix \eqref{AA}. Put $A_{ij}=\omega^{a_{ij}}$, and $H=\overline{A}\circ B$, $H_{ij}=\omega^{-a_{ij}}\cdot\zeta^{b_{ij}}$. Then, we have $B=J$.
\end{proposition}
\begin{proof}
Consider the linear system of $376$ equations \eqref{eq1}--\eqref{eqlin} arising in Lemma~\ref{lem4}. Since $B$ is in normal form, $b_{1j}=0$ for every $j\in\{1,2,\dots,16\}$. The coefficient matrix of this system in the $240$ remaining variables $b_{ij}$, $i\in\{2,3,\dots,16\}$, $j\in\{1,2,\dots,16\}$ can be reconstructed explicitly using the matrix $A$ from \eqref{AA}. Then, performing Gaussian elimination over $\mathrm{GF}(5)$ (by a conveniently chosen computer algebra system, say Mathematica\footnote{Wolfram Research, Inc., Mathematica, Version 14.3, Champaign, IL (2025).}), one may verify that the rank of this matrix is $240$. In particular, all the variables must be simultaneously $0$.
\end{proof}
\begin{theorem}\label{T29}
There do not exist $\mathrm{BH}(16,15)$ matrices.
\end{theorem}
\begin{proof}
Using the notations of the statement of Proposition~\ref{prop2}, let $H=\overline{A}\circ B$ be such a hypothetical matrix in normal form. Then, $H=\overline{A}\circ J=\overline{A}$ is a $\mathrm{BH}(16,3)$ matrix which is impossible since $3\nmid 16$.
\end{proof}
\begin{remark}
We are confident that similar nonexistence results can be obtained for other suitably chosen pairs $(n,q)$ as well by simply following our approach. However, as $n$ grows, it becomes increasingly difficult to analyse (and cumbersome to document) the arising projection matrices (cf.~\eqref{AA}), and it seems that new ideas are required for further significant progress.
\end{remark}
\section{The nonexistence of BH(17,15) matrices}\label{sec3}
The aim of this section is to prove Theorem~\ref{thm:ferenc}, and consequently the nonexistence of $\mathrm{BH}(17,15)$ matrices. We will show that carefully chosen parameters would give rise to $\mathrm{BH}(n,3p)$ matrices $H$ having a highly structured projection $H^{\circ p}$, whose existence in turn would lead to a strongly regular graph \cite{srgbook}. Using this correspondence, known nonexistence results of these graphs imply the nonexistence of Butson matrices. We begin by analysing the Gram structure of the projections.
\begin{lemma}\label{prop31}
Let $p\equiv 2\ (\mathrm{mod}\ 3)$ be a prime number, let $s\in\{2,3,\dots,p-1\}$ be a parameter, and write $n=p+3s$. Let $H$ be a $\mathrm{BH}(n,3p)$ matrix. Then, there exists a cube root matrix $A$ equivalent to $H^{\circ p}$, such that for every $i,j\in\{1,2,\dots,n\}$, we have
\begin{align}
\prod_{k=1}^nA_{ik}=\prod_{k=1}^nA_{kj}=1,\label{prod1}\\
AA^\ast=A^\ast A=3sI+pJ,\label{prod2}\\
AJ=JA=\alpha J,\quad |\alpha|^2=3s+np\label{prod3}.
\end{align}
\end{lemma}
\begin{proof}
Let $i<j\in\{1,2,\dots,n\}$. Since $n<4p$, the inner product $\left\langle H_i,H_j\right\rangle$ decomposes as a vanishing sum of exactly one $p$-cycle and $s$ sets of $3$-cycles \cite{LamLeung} similarly as seen in the proof of Lemma~\ref{l2}. Therefore every quotient row $(H_i\circ \overline{H_j})^{\circ p}=A_i\circ\overline{A_j}$ contains $s$ entries of $1$, $\omega$, and $\omega^2$ each, and in addition, $p$ more entries equal to $\omega^{r_{ij}}\in\Omega_3$ (exponents are understood modulo $3$). Observe that we maintain this structure while replacing $A$ with another one equivalent to it.

Next we show that the matrix $A$ is equivalent to one in which the product of the entries in every row and every column equals $1$. Let us denote the product of the entries in the $i$th row of $A$ by $x_i:=\prod_{k=1}^n A_{ik}\in\Omega_3$. Multiplying this row by $x_i$ results in an equivalent cube root matrix where this row product becomes $x_i^{n+1}=x_i^{p+1+3s}=1$. After repeating this procedure with every row, consider the product of the entries in the $j$th column ($j\in\{1,2,\dots,n\}$) of the resulting matrix $A'$, denoted by $y_j:=\prod_{k=1}^nA'_{kj}\in\Omega_3$. Since $\prod_{j=1}^n y_j=\prod_{i=1}^n x_i^{n+1}=1$, we may multiply every column of $A'$ by their respective product $y_j$, while keeping the product of the entries in each row unchanged. This proves \eqref{prod1}.

Therefore, we may assume without loss of generality, that the product of the entries in every row and column of $A$ equals $1$. We may compare this with the product of the entries in a quotient row $A_i\circ \overline{A_j}$ ($i<j\in\{1,2,\dots n\}$), as
\[1=\prod_{k=1}^n A_{ik}\cdot \prod_{k=1}^n \overline{A_{jk}}=\prod_{k=1}^n (A_{i}\circ \overline{A_{j}})_k=1^s\cdot \omega^s\cdot\omega^{2s}\cdot\omega^{pr_{ij}}=\omega^{pr_{ij}+3s}=\omega^{2r_{ij}}.\]
Therefore $r_{ij}=0$ for every $i<j\in\{1,2,\dots,n\}$ and consequently $\left\langle A_i,A_j\right\rangle = p\omega^{r_{ij}}=p$ as claimed. Since the same argument applies to the quotient rows of $A^T$, we have proved \eqref{prod2}.

Multiplying either side of \eqref{prod2} by $A$, we obtain $AA^\ast A=3sA+pJA=3sA+pAJ$, hence $AJ=JA=\alpha J$, where $\alpha\in\mathbb{Q}(\omega)$ is the common row and column sum of $A$. Then, using once again \eqref{prod2}, we have $n|\alpha|^2J=JA(JA)^\ast=J(3sI+pJ)J=n(3s+np)J$, and we have \eqref{prod3}.
\end{proof}
\begin{proposition}\label{newprop}
Let $p\equiv 2\ (\mathrm{mod}\ 3)$ be a prime number, let $s$ be an integer parameter such that
\begin{equation}\label{sdef}
\sqrt{4p+5}\leq 2s-3\leq 2p-5,
\end{equation}
and write $n=p+3s$. Let $H$ be a $\mathrm{BH}(n,3p)$ matrix and assume that
\begin{equation}\label{kdef}
k=\frac{n-1-\sqrt{(n-1)(n-4s+3)}}{2}
\end{equation}
is an integer. If the system of equations
\begin{equation}\label{syseqs}\begin{cases}
a+b+c=n\\
ab+ac+bc=s(n-1)
\end{cases}
\end{equation}
admits the only integer solution $(a,b,c)=(n-k-1,k,1)$ with $a\geq b\geq c\geq 0$, then there exists a strongly regular graph with parameters $(n,k,k-s,k-s+1)$.
\end{proposition}
\begin{proof}
On the one hand, a simple calculation shows that the lower bound on $s$ in \eqref{sdef} implies that $k\geq s$, hence the parameters of the claimed graph are well-defined.

On the other hand, $s\leq p-1$, and therefore we may apply Lemma~\ref{prop31} to obtain the matrix $A$ equivalent to $H^{\circ p}$. By \eqref{prod3}, we have $AJ=JA=\alpha J$, where $\alpha\in\mathbb{Q}(\omega)$ is the common row and column sum of $A$. Assume that the triple $(a,b,c)\in\mathbb{Z}^3$ counts the number of entries within the first row of $A$ equal to $1$, $\omega$, and $\omega^2$, respectively. Then $\alpha=a+b\omega+c\omega^2$, where $a+b+c=n$. Calculating $|\alpha|^2$ in two different ways (cf.~\eqref{prod3}) gives
\[|\alpha|^2=n^2-3(ab+ac+bc)=3s+np,\]
which reduces to the system \eqref{syseqs}.

Since the row and column sums are equal, the triple $(a,b,c)$ is independent of the chosen row of $A$ (or that of $A^T$). Indeed, assuming that $\alpha=a+b\omega+c\omega^2=a'+b'\omega+c'\omega^2$ and using that $\omega^2=-1-\omega$, we have $(a',b',c')=(a+z,b+z,c+z)$ for some $z\in\mathbb{Z}$. However, since $a+b+c=a'+b'+c'=n$, $z=0$ and the two triples coincide. We may assume that $a=n-k-1$ by multiplying $A$ by a cube root if it is necessary. Furthermore, conjugating $A$ interchanges the role of $b$ and $c$, hence we may have $(a,b,c)=(n-k-1,k,1)$ in every row and column. These operations leave \eqref{prod2} and \eqref{prod3} invariant.

The positions of the unique $\omega^2$ in each row and column give rise to a permutation matrix. Therefore, we may assume without loss of generality, that the columns of the matrix $H$ (and hence, that of $A$), have been permuted in advance so that the entries $\omega^2$ are on the main diagonal. These operations once again leave \eqref{prod2} and \eqref{prod3} invariant. Let $X$ denote the $\{0,1\}$-matrix indicating the positions of $\omega$ in $A$. Then, we have
\[
A=J+(\omega-1)X+(\omega^2-1)I,
\qquad X e=X^T e=k e.
\]
Comparing the expansion of $AA^\ast$ with \eqref{prod2} and collecting the coefficients of $\omega$ and $1$ yields $X=X^T$ and $X^2=(s-1)I-X+(k-s+1)J$, respectively. In particular, $X$ is the adjacency matrix of a $k$-regular graph. Let $f\in\mathbb{R}^n$ be an eigenvector of $X$ with eigenvalue $\theta$, orthogonal to $e$. Then $\theta^2+\theta-(s-1)=0$. Since the discriminant of this quadratic equation is $1+4(s-1)\neq 0$, $X$ has exactly two restricted eigenvalues \cite{srgbook}, and hence it is the adjacency matrix of a strongly regular graph.
\end{proof}
\begin{theorem}\label{thm:ferenc}
Let $p\equiv 2\ (\mathrm{mod}\ 3)$ be a prime number, let $s$ be an integer parameter satisfying \eqref{sdef}, and write $n=p+3s$. Assume that $k$ as in \eqref{kdef} is an integer, and $4s-3$ is not the square of an integer. If the system of equations \eqref{syseqs} admits the only integer solution $(a,b,c)=(n-k-1,k,1)$ with $a\geq b\geq c\geq 0$, then there do not exist $\mathrm{BH}(n,3p)$ matrices.
\end{theorem}
\begin{proof}
Suppose to the contrary that a $\mathrm{BH}(n,3p)$ matrix exists, and let $X$ be the adjacency matrix of the strongly regular graph with parameters $(n,k,k-s,k-s+1)$, coming from Proposition~\ref{newprop}. The eigenvalues of $X$ are $k$ (of multiplicity one), and the roots of $\theta^2+\theta-(s-1)=0$, which read $\theta_{1,2}=(-1\pm\sqrt{4s-3})/2$. Since $4s-3$ is not a perfect square, $\operatorname{tr}(X)=0$ forces the multiplicities of $\theta_1$ and $\theta_2$ to be equal to $(n-1)/2$. This is impossible, since $n-4s+3\geq 4$, and consequently
\[0=\operatorname{tr}(X)=k-(n-1)/2=-\sqrt{(n-1)(n-4s+3)}/2\leq -\sqrt{n-1}.\qedhere\]
\end{proof}
\begin{corollary}\label{cor}
There do not exist $\mathrm{BH}(17,15)$ matrices.
\end{corollary}
\begin{proof}
Apply Theorem~\ref{thm:ferenc} with $p=5$ and $s=4$.
\end{proof}
\begin{remark}
Proposition~\ref{newprop} has further applications. For example, a $\mathrm{BH}(50,3\cdot 29)$ matrix would lead to a cube root projection matrix $A$ which corresponds to an existing strongly regular graph with parameters $(50,7,0,1)$: the Hoffman--Singleton graph \cite{srgbook}. Unfortunately, a very similar lifting argument with minor modifications as described in Lemma~\ref{lem4} shows that no such Butson matrix can exist. It would be interesting to see a case where this lifting works and leads to a new Butson matrix.
\end{remark}
\section*{Acknowledgements and Statement of AI use}
During early stages of this research, OpenAI ChatGPT (GPT-5.6 Pro and GPT-6 Astra) were used by some of the authors. However, this entire manuscript has been rewritten from scratch, and all proofs have been independently rederived. The authors take full responsibility for all claims made in this manuscript.

\textsc{D. C.} was supported by the Ministry
of Innovation and Technology NRDI Office within the framework of the Artificial Intelligence National Laboratory (RRF-2.3.1-21-2022-00004), and NKFIH-154121. \textsc{I. N.} was supported by the grant NKFIH-154121. \textsc{F. Sz.} was supported in part by JSPS KAKENHI Grant Number 24K06829.

\end{document}